\documentclass[a4paper]{article}
\usepackage[latin1]{inputenc}
\usepackage{amscd}
\usepackage[english]{babel}
\usepackage{amsmath,amssymb,amsfonts}
\usepackage{syntonly}

\usepackage{color}
\usepackage[outerbars,color]{changebar}
\usepackage{ifthen}

\newboolean{usecolor}
\setboolean{usecolor}{true}

\ifthenelse{\boolean{usecolor}}
{
  \definecolor{colore}{cmyk}{0,1,0.6,0}
  \definecolor{coloregen}{cmyk}{0.7,0,1,0}
  \definecolor{coloresimo}{cmyk}{1,0.6,0,0}
}
{
  \definecolor{colore}{cmyk}{0,0,0,1}
  \definecolor{coloregen}{cmyk}{0,0,0,1}
  \definecolor{coloresimo}{cmyk}{0,0,0,1}
}

\DeclareFixedFont{\fiverm}{OT1}{cmr}{m}{n}{5}

\newcommand{\Z}[0]{\mathbb{Z}}
\newcommand{\R}[0]{\mathbb{R}}
\newcommand{\C}[0]{\mathbb{C}}

\newcommand{\Cal}[1]{\mathcal{#1}}

\newcommand{\eq}[1][r]
       {\ar@<-3pt>@{-}[#1]
        \ar@<-1pt>@{}[#1]|<{}="gauche"
        \ar@<+0pt>@{}[#1]|-{}="milieu"
        \ar@<+1pt>@{}[#1]|>{}="droite"
        \ar@/^2pt/@{-}"gauche";"milieu"
        \ar@/_2pt/@{-}"milieu";"droite"}
\newcommand{\imm}[1][r] {\ar@{^{(}->}[#1]}
\renewcommand{\ni}[0]{\noindent}

\newtheorem{df}{Definition}[section]
\newtheorem{teo}{Theorem}
\newtheorem{prop}[df]{Proposition}
\newtheorem{lem}[df]{Lemma}

\newtheorem{conj}[df]{Conjecture}
\newtheorem{rmk}[df]{Remark}

\newcommand{\wdt}{\widetilde}

\newcommand{\rx}{\mathcal{R}_X}

\input prepictex
\input pictex
\input postpictex

\title{The homotopy type of toric arrangements}
\author {L. Moci \and  \and S. Settepanella}
\author{Luca {\sc Moci}\footnote{Institut fuer Mathematik, Technische Universitaet Berlin,
Strasse des 17. Juni, 136, Berlin. moci@math.tu-berlin.de (This work was partially supported by a Sofia Kovalevskaya Research Prize of
Alexander von Humboldt Foundation awarded to Olga Holtz.)}
\and Simona {\sc Settepanella}\footnote{LEM, Scuola Superiore Sant'Anna, Pisa, Italy. s.settepanella@sssup.it (Thanks to financial support from the European Commission 6th FP (Contract CIT3-CT-2005-513396), Project: DIME - Dynamics of Institutions and Markets in Europe) }}

\begin{document}

\maketitle
\begin{abstract}
A toric arrangement is a finite set of hypersurfaces in a complex torus, every hypersurface being the kernel of a character. In the present paper we build a CW-complex homotopy equivalent to the arrangement complement, with a combinatorial description similar to that of the well-known Salvetti complex. If the toric arrangement is defined by a Weyl group, we also provide an algebraic description, very handy for cohomology computations. In the last part we give a description in terms of tableaux for a toric arrangement appearing in robotics.
\end{abstract}

\begin{center}
\end{center}

\begin{center}
\end{center}

\section*{Introduction}
A toric arrangement is a finite set of hypersurfaces in a complex torus $T=(\C^*)^n$, in which every hypersurface is the kernel of a character $\chi \in X \subset Hom(T,\C^*) $ of $T$.
Let $\rx$ be the complement of the arrangement: its geometry and topology have been studied by many authors, see for instance  \cite{L2}, \cite{li}, \cite{DP}, \cite{Mo2}.
In particular, in \cite{Lo} and \cite{DP} the De Rham cohomology of $\rx$ has been computed, and recently in \cite{Mo3} a \emph{wonderful model} has been built.

In the present paper we build a topological model $\Cal{S}$ for $\rx$. This
model is a regular CW-complex, similar to the one introduced by Salvetti (\cite{Sa1}) for the complement of hyperplane arrangements.
Indeed for a wide class of arrangements, which we call \textit{thick},
the cells of $\Cal{S}$ are in bijection with  pairs $[C \prec F]$, where
$C$ is a \textit{chamber} of the \textit{real} toric arrangement and
$F$ is a \textit{facet} adjacent to it (according to the definitions
given in Section \ref{DefTo}).

\bigskip

The model $\Cal S$ is well suited for homology and homotopy computations, which will be developed in future papers (see for instance \cite{simo2}).
Furthermore, 
the jumping loci in the local system cohomology of a CW-complex are affine algebraic varieties. In the theory of hyperplane arrangements such objects, called \emph{characteristic varieties}, proved to be of fundamental importance. It is then a remarkable fact that the characteristic varieties can be defined also in the toric case.

In Section 2 we focus on a toric arrangement associated to an affine
Weyl group $\widetilde{W}$.
In this case the chambers are in bijection with the elements
of the corresponding finite Weyl group W, and the cells of $\Cal S$ are
given by the pairs $(w, \Gamma)$, where $w \in W$ and $\Gamma$ is a
proper subset of the set $S$ of generators of  $\widetilde{W}$.
This generalizes a construction introduced in \cite{Sa} and \cite{boss3}.

In the last Section we give a description of the facets of the real toric arrangement defined by the Weyl group $\wdt A_n$ in the torus corresponding to the root lattice. This description in terms of Young tableaux turns out to be interesting since it coincides with the complex
describing the space of all periodic legged gaits of a robot body (see \cite{KGCohen}).

\bigskip

\paragraph{Acknowledgements}
We are grateful to the organizers of the research program "Configuration Spaces: Geometry, Combinatorics and Topology" at Centro De Giorgi (Pisa), which provided us a significant occasion to work together. In particular we wish to thank Fred Cohen and Mario Salvetti for several valuable suggestions. We also thank Priyavrat Deshpande for many stimulating conversations we had while we were completing the present paper.

\section{The CW-complex}

\subsection{Main definitions}
Let $T=(\C^*)^n$ be a complex torus and $X \subset Hom(T,\C^*) $ be a finite set of
characters of $T$. The kernel of every $\chi \in X$ is a hypersurface
of $T$:
\begin{equation*}
H_{\chi}:=\{t \in T \, \mid \, \chi(t)=1\}.
\end{equation*}
Then $X$ defines on T the \textit{toric arrangement}:
\begin{equation*}
\mathcal{T}_X:=\{ H_{\chi} , \chi \in X\}.
\end{equation*}
Let $\rx$ be the \textit{complement } of the arrangement:
\begin{equation*}
\rx:= T \setminus \bigcup_{\chi \in X} H_{\chi}.
\end{equation*}
Let $\pi: V \longrightarrow T$ be the universal covering of $T$. Then
$V$ is a complex vector space of rank $n$, and $\pi$ is the quotient
map $\pi: V \longrightarrow V/ \Lambda$, where $\Lambda$ is a lattice
in $V$. Then the preimage $\pi^{-1}(H_{\chi})$ of a hypersurface
$H_{\chi} \in \mathcal{T}_X$ is the union of an infinite family of parallel hyperplanes. Thus
$$ \Cal A_X:=\{ H \mbox{ hyperplane of } V \mid \exists \chi \in X
\mbox{ s.t. } \pi(H)=H_{\chi}\}$$
is a periodic affine hyperplane arrangement in V. Let $\mathcal{M}_X$ be its
complement:
\begin{equation*}
\mathcal{M}_X:= V \setminus \bigcup_{\chi \in X}\pi^{-1}(H_{\chi}).
\end{equation*}
By definition, $\pi$ maps $\mathcal{M}_X$ on $\rx$. Moreover the equations
defining the hyperplanes in $\Cal A_X$ can always be assumed to have
integral (hence real) coefficients since they are given by elements of
$\Lambda$. Thus by \cite{Sa1} there is an (infinite) CW-complex
$\widetilde{\Cal S} \subset \mathcal{M}_X$ and a map $\varphi: \mathcal{M}_X
\longrightarrow \widetilde{\Cal S}$ giving a homotopy equivalence. \\
Furthermore, we can build $\widetilde{\Cal S}$ in such a way that it is
invariant under the action of translation in $\Lambda$: for instance
by building the cells relative to a fundamental domain and then inductively,
defining for each cell above the other cells of its $\Lambda$-orbit by translation.
Thus $\pi(\widetilde{\Cal S})$ is a finite CW-complex, which will be
denoted by $\mathcal{S}$, and the image of every cell of $\widetilde{\Cal S}$
is a cell of $\mathcal{S}$. Moreover, since $\varphi$ is $\Lambda-$equivariant,
it is well defined the map
\begin{equation*}
\varphi_{\pi}(t):= (\pi \circ \varphi)(\pi^{-1}(t))
\end{equation*}
which makes the following diagram commutative:
\begin{equation} \begin{array}{ccc}\label{complesso}
\mathcal{M}_X &\xrightarrow{\varphi} &\widetilde{\Cal S}  \\
\pi \downarrow & &\pi \downarrow \\
\rx & \xrightarrow{\varphi_{\pi}} & \Cal S
\end{array}
\end{equation}

\begin{lem} The map $\varphi_{\pi}$ is a homotopy equivalence between
  $\rx$ and $\Cal S$.
 \end{lem}

\textbf{Proof.} The map $\varphi$ is a homotopy equivalence hence, by
definition, there is a continuous map $\psi: \widetilde{\Cal S}
\rightarrow \mathcal{M}_X$
such that $\psi \varphi$ is homotopic to the
identity map $id_{\mathcal{M}_X}$ and $\varphi \psi$ is homotopic to
$id_{\widetilde{\Cal S}}$. Namely, since $\widetilde{\Cal S}$ is a deformation retract, the homotopy inverse $\psi$ is simply the inclusion map, which is clearly $\Lambda-$equivariant. Hence the map
$$\psi_\pi(t):= (\pi \circ \psi)(\pi^{-1}(t))$$
is well defined and makes the following diagram commutative:
\begin{equation} \begin{array}{ccc}
\widetilde{\Cal S}  &\xrightarrow{\psi} & \mathcal{M}_X\\
\pi \downarrow & &\pi \downarrow \\
\Cal S & \xrightarrow{\psi_{\pi}}  & \rx .
\end{array}
\end{equation}

Let $I=[0,1]$ be the unit interval and $F: \mathcal{M}_X \times I \rightarrow
\mathcal{M}_X$ be the continuous map such that $F(x,0) = \psi(\varphi(x))$ and
$F(x,1)=id_{M_{X}}(x)$. Again, since $F$ is $\Lambda-$equivariant, we can define the map:
$$F_\pi(t):= (\pi \circ F)(\pi^{-1}(t))$$
In this way we get the commutative diagram:
\begin{equation} \begin{array}{ccc}
\mathcal{M}_X \times I&\xrightarrow{F} & \mathcal{M}_X  \\
\pi \downarrow & &\pi \downarrow \\
\rx \times I& \xrightarrow{F_{\pi}} & \rx.
\end{array}
\end{equation}

By construction map $F_{\pi}$ is a continuous map such that $F_{\pi}(x,1)=id_{\rx}$ and
$$F_{\pi}(x,0)= (\psi\varphi)_{\pi}(x)=\pi\psi \varphi\pi^{-1}(x)=\pi\psi \pi^{-1}\pi\varphi\pi^{-1}(x)=\psi_\pi\circ\varphi_\pi(x).$$

Hence $F_{\pi}$ gives the required homotopy equivalence. $\qquad \square$

\subsection{Salvetti complex for affine arrangements}
In order to describe the structure of $\Cal S$, we now have to focus on the real counterparts of the complex arrangements above.\\
Let $V_{\R}$ be the real part of $V$. In other words, let $V_{\R}\doteq \R^n$ be a real vector space, and let $V\doteq V_{\R}\otimes_\R \C$ be its complexification. Then we identify $V_{\R}$ with a subspace of $V$ via the map $v\mapsto v\otimes 1$.\\
Let $\Cal{A}_{X,\R}$ be the
corresponding hyperplane arrangement on $V_{\R}$ and $\mathcal{M}_{X,\R}=\mathcal{M}_{X}\cap V_{\R}$ its
complement.
Since the image of $\R$ under the map $\C \longrightarrow \C / \Z
\xrightarrow{\sim} \C^*$ is the circle
$$S^1:=\{z \in \C \, \mid \, \mid z \mid =1\}$$
we have that the image of $V_{\R}$ under the map $\pi:V \rightarrow V/\Lambda
\xrightarrow{\sim} T$ is a compact torus $T_{\R} \subset T$. A \emph{real toric
arrangement} $\Cal T_{X,\R}$ is naturally defined on $T_{\R}$ with
hypersurfaces $H_{\chi,\R}:=H_{\chi} \cap T_{\R}$ and complement
$\mathcal{R}_{X,\R} = \rx \cap T_{\R}$. Furthermore $\pi$ restricts to universal
covering map $\pi: V_{\R} \longrightarrow T_{\R}$ and
$\pi(\mathcal{M}_{X,\R})=\mathcal{R}_{X,\R}$.

\bigskip

We recall the following definitions:
\begin{enumerate}
\item a \textit{chamber} of $\Cal A_{X,\R}$ is a connected component
  of $\mathcal{M}_{X,\R}$;
\item a \textit{space} of $\Cal A_{X,\R}$ is an intersection of
  elements in $\Cal A_{X,\R}$;
\item a \textit{facet} of $\Cal A_{X,\R}$ is the intersection of a space
  and the closure of a chamber.
\end{enumerate}
Let $\mathbf{S}:=\{\widetilde{F}^k\}$ be the stratification of $V_{\mathbb{R}}$ into facets $\widetilde{F}^k$
induced by the arrangement $\Cal{A}_{X,\R}$ (see \cite{Bou}), where superscript $k$ stands
for codimension.

Then the $k$-cells of the complex $\widetilde{\Cal S}$ described in
\cite{Sa1}  bijectively correspond to pairs
$$[\widetilde{C} \prec \widetilde{F}^k]$$
where $\widetilde{C}=\widetilde{F}^0$ is a chamber of $\mathbf S$ and
$\widetilde{F}^i \ \prec \widetilde{F}^j  \Leftrightarrow
clos(\widetilde{F}^i)\supset \widetilde{F}^j$ is the standard partial
ordering in  $\mathbf S$ (see also \cite{OT}).

Let $|\widetilde{F}|$ be the affine subspace spanned by $\widetilde{F},$ and let us consider
the subarrangement
$$\Cal A_{\widetilde{F}}\ =\ \{ H \in \Cal{A}_{X,\R} \ :\ \widetilde{F}\subset H\}.$$
A cell $[\widetilde{C}\prec \widetilde{F}^k]$ is in the boundary of $[\widetilde{D}\prec \widetilde{G}^j]$ ($k<
j$) if and only if

\begin{equation}\label{bordoaffine}
\begin{split}
& \mbox{i) } \widetilde{F}^k\prec \widetilde{G}^j \\
& \mbox{ii) } \mbox{ the chambers } \widetilde{C} \mbox{ and } \widetilde{D} \mbox{ are contained in the same chamber of } \Cal A_{\widetilde{F}^k}.
\end{split}
\end{equation}

 Previous conditions are equivalent to say that
  $\widetilde{C}$ is the chamber of $\Cal{A}_{X,\R}$ which is the "closest" to $\widetilde{D}$
among those which contain $\widetilde{F}^k$ in their closure. The standard notation
$[\widetilde{C}\prec \widetilde{F}^k] \in \partial_{\widetilde{\Cal S}}[\widetilde{D}\prec \widetilde{G}^j]$ will be used.

\bigskip

 {It is a simple remark that the above description of the Salvetti complex 
$\widetilde{\Cal S}$ is $\Lambda$-invariant. Indeed each translation $t \in \Lambda$ acts on the stratification $\mathbf{S}:=\{\widetilde{F}^k\}$ sending a $k$-facet $F^k$ into the $k$-facet 
$t.F^k$. Then the translation $t$ acts on $\widetilde{\Cal S}$ sanding a $k$-cell $[C \prec F^k]$
in the $k$-cell $[t.C \prec t.F^k]$.}

\subsection{Salvetti Complex for toric arrangements}\label{DefTo}

In order to give a similar description for $\Cal S$, we introduce the
following definitions:
\begin{enumerate}
\item a \textit{chamber} of $\Cal T_{X,\R}$ is a connected component
  of $R_{X,\R}$;
\item a \textit{layer} of $\Cal T_{X,\R}$ is a connected component of
  an intersection of elements of $\Cal T_{X,\R}$;
\item a \textit{facet} of $\Cal T_{X,\R}$ is an intersection of a layer
  and the closure of a chamber.
\end{enumerate}

\begin{lem}\label{camereefacce}

~

\begin{enumerate}
\item If $\widetilde{C}$ is a chamber of $\Cal A_{X,\R}$,
 $\pi(\widetilde{C})$ is a chamber of $\Cal T_{X,\R}$;
\item If $\widetilde{L}$ is a space of $\Cal A_{X,\R}$,
 $\pi(\widetilde{L})$ is a layer of $\Cal T_{X,\R}$;
\item If $\widetilde{F}$ is a facet of $\Cal A_{X,\R}$,
 $\pi(\widetilde{F})$ is a facet of $\Cal T_{X,\R}$;
\end{enumerate}
\end{lem}

\textbf{Proof.} The first statement is clear, as well as the second
one since $\pi(\widetilde{L})$ must be connected. The third claim is a
direct consequence of the previous two.~$~\square~$

\bigskip


Now, let us consider pairs $$[C \prec F^k]$$
where $C=F^0$ is a chamber of $\Cal T_{X,\R}$,  $F^k$ a
$k$-codimensional facet of
$\Cal T_{X,\R}$ and $F^i \prec F^j  \Leftrightarrow
clos(F^i)\supset F^j$.

By Lemma \ref{camereefacce} the quotient map $\pi(\wdt F)$ of a facet
is still a facet in the real torus and, because of the surjectivity of $\pi$, we get that any facet $F$ in $\Cal{T}_{X, \R}$ is the image
$F=\pi(\wdt F)$ of an affine one.

We notice that
$$\pi([\wdt C \prec \wdt F])=\pi([\wdt D \prec \wdt G])  \Longrightarrow [\pi(\wdt C) \prec \pi(\wdt F)] =[\pi(\wdt D) \prec \pi(\wdt G)]. $$
Indeed if $\pi([\wdt C \prec \wdt F])=\pi([\wdt D \prec \wdt G])$ there is a translation
$t \in \Lambda$ which sends $[\wdt C \prec \wdt F]$ in $[\wdt D \prec \wdt G]$. As a simple
consequence $\wdt D=t. \wdt C$ and $\wdt F =t.\wdt G$, i.e. $\pi(\wdt C)=\pi(\wdt D)$ and
$\pi(\wdt F)=\pi(\wdt D)$.

Then there is a natural surjective map from the cells of $\Cal S$ to the set of pairs $[C \prec F]$, but this map in general is not injective.
Let us consider the simple example defined by $\Cal A=\{x \in \R \mid x \in \Z \}$.

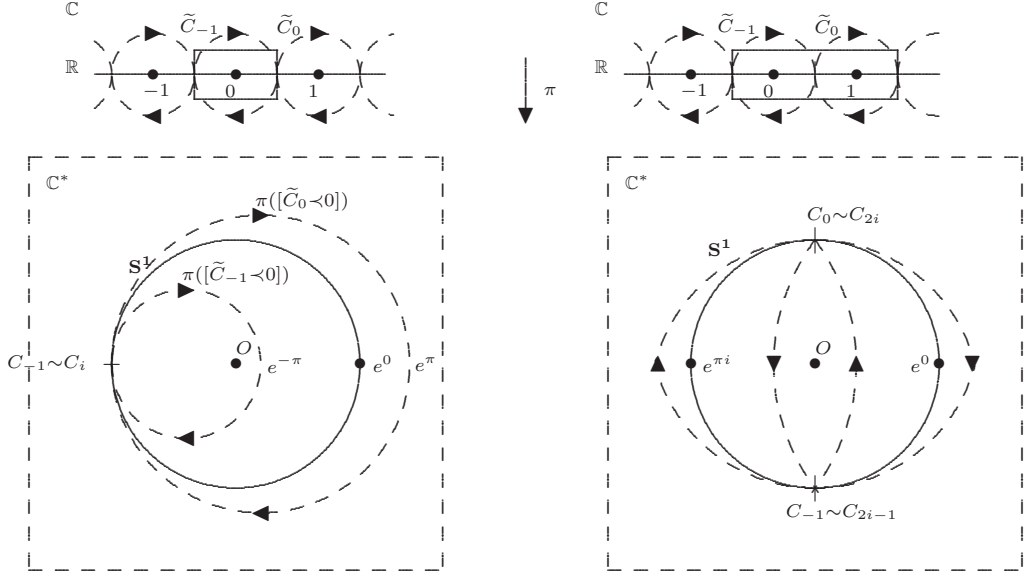
\begin{equation*}
\beginpicture
\setcoordinatesystem units <1.1cm,1.1cm>         
\setplotarea x from -6 to 5, y from -1 to 6    
\put{$\bullet$} at  -4.5  5
\put{$\bullet$} at  -3.5  5
\put{$\bullet$} at  -2.5   5
\put{$\bullet$} at  2  5
\put{$\bullet$} at  3  5
\put{$\bullet$} at  4  5
\put{$\mid$} at  -4  5
\put{$\mid$} at  -3  5
\put{$\mid$} at  2.5  5
\put{$\mid$} at  3.5  5
\put{$\mid$} at  4.5  5
\put{$\blacktriangleright$} at -4.5 5.5
\put{$\blacktriangleright$} at -3.5 5.5
\put{$\blacktriangleright$} at -2.5 5.5
\put{$\blacktriangleleft$} at  -4.5  4.5
\put{$\blacktriangleleft$} at  -3.5 4.5
\put{$\blacktriangleleft$} at  -2.5 4.5
\put{$\blacktriangleright$} at 2 5.5
\put{$\blacktriangleright$} at 3 5.5
\put{$\blacktriangleright$} at 4 5.5
\put{$\blacktriangleleft$} at  2 4.5
\put{$\blacktriangleleft$} at  3 4.5
\put{$\blacktriangleleft$} at  4 4.5
\put{$\scriptstyle{-1}$}[r] at -4.3 4.8
\put{$\scriptstyle{0}$}[r] at -3.5 4.8
\put{$\scriptstyle{1}$}[r] at -2.5 4.8
\put{$\scriptstyle{-1}$}[r] at 2.2 4.8
\put{$\scriptstyle{0}$}[r] at 3 4.8
\put{$\scriptstyle{1}$}[r] at 4 4.8
\put{$\scriptstyle{\wdt C_{-1}}$}[r] at -3.7 5.6
\put{$\scriptstyle{\wdt C_{0}}$}[r] at -2.7 5.6
\put{$\scriptstyle{\mathbf{\C^*}}$}[r] at -5.5 3.7
\put{$\scriptstyle{\mathbf{S^1}}$}[r] at -4.5 2.7
\put{$\scriptstyle{\wdt C_{-1}}$}[r] at 2.8 5.6
\put{$\scriptstyle{\wdt C_{0}}$}[r] at 3.8 5.6
\put{$\scriptstyle{\mathbf{\C^*}}$}[r] at 1.5 3.7
\put{$\scriptstyle{\mathbf{S^1}}$}[r] at 2.5 2.9
\plot  -5.2  5   -1.7  5 /
\put{$\scriptstyle{\mathbf{\R}}$}[r] at -5.4 5.1
\put{$\scriptstyle{\mathbf{\C}}$}[r] at -5.4 5.8
\plot  1.2  5   5  5 /
\put{$\scriptstyle{\mathbf{\R}}$}[r] at 1 5.1
\put{$\scriptstyle{\mathbf{\C}}$}[r] at 1 5.8
\put{$\bullet$} at  -3.5  1.5
\put{$\scriptstyle{O}$}[r] at -3.3 1.7
\put{$\bullet$} at  -2  1.5
\put{$\scriptstyle{e^0}$}[r] at -1.6 1.5
\put{$\bullet$} at  3.5  1.5
\put{$\scriptstyle{O}$}[r] at 3.7 1.7
\put{$\bullet$} at  2  1.5
\put{$\scriptstyle{e^{\pi i}}$}[r] at 2.5 1.5
\put{$\bullet$} at  5  1.5
\put{$\scriptstyle{e^0}$}[r] at 4.9 1.5
\put{$-$} at  -5  1.5
\put{$\scriptstyle{C_{-1} \sim C_{i}}$}[r] at -5.3 1.5
%
\put{$\mid$} at  3.5  3
\put{$\scriptstyle{C_{0} \sim C_{2i}}$}[r] at 4.3 3.3
\put{$\mid$} at  3.5  0
\put{$\scriptstyle{C_{-1} \sim C_{2i-1}}$}[r] at 4.5 -0.3
\circulararc 360 degrees from -5  1.5  center at -3.5 1.5
\circulararc 360 degrees from 2  1.5  center at 3.5  1.5
\plot  -4 5.3   -3 5.3 /
\plot  -4 4.7    -3 4.7 /
\plot  -4 5.3   -4 4.7 /
\plot  -3 5.3  -3  4.7 /
\plot  2.5 5.3   4.5 5.3 /
\plot  2.5 4.7    4.5 4.7 /
\plot  2.5 5.3   2.5 4.7 /
\plot  4.5 5.3  4.5  4.7 /
\plot  0 5.2   0 4.5 /
\put{$\blacktriangledown$} at 0 4.5
\put{$\scriptstyle{\pi}$} at 0.3 4.8
\setdashes
\circulararc 360 degrees from -4  5  center at -4.5 5
\circulararc 360 degrees from -3  5  center at  -3.5 5
\circulararc 360 degrees from -3  5  center at  -2.5 5
\circulararc 60 degrees from -5  5  center at  -5.5 5
\circulararc -60 degrees from -5  5  center at  -5.5 5
\circulararc 80 degrees from -2  5  center at  -1.5 5
\circulararc -80 degrees from -2  5  center at  -1.5 5
\circulararc 360 degrees from 2.5  5  center at 2 5
\circulararc 360 degrees from 2.5  5  center at 3  5
\circulararc 360 degrees from 3.5  5  center at 4  5
\circulararc 90 degrees from 4.5  5  center at 5  5
 \circulararc -90 degrees from 4.5  5  center at 5  5
 \circulararc 50 degrees from 1.5  5  center at 1  5
 \circulararc -50 degrees from 1.5  5  center at 1  5
%
\circulararc 360 degrees from -5  1.5  center at  -3.2 1.5
\put{$\scriptstyle{e^{\pi}}$} at -1.2 1.5
\put{$\blacktriangleright$} at -3.2 3.3
\put{$\blacktriangleleft$} at -3.2 -0.3
\circulararc 360 degrees from -5  1.5  center at -4.1 1.5
\put{$\scriptstyle{\pi([\wdt C_{-1} \prec 0])}$} at -3.5 2.6
\put{$\scriptstyle{e^{-\pi}}$} at -2.9 1.5
\put{$\blacktriangleright$} at -4.1 2.4
\put{$\blacktriangleleft$} at -4.1 0.6
\put{$\scriptstyle{\pi([\wdt C_{0} \prec 0])}$} at -2.7 3.5
\circulararc -75 degrees from 3.5  3  center at 3.5 1
\circulararc 75 degrees from 3.5  0  center at 3.5 2
\put{$\blacktriangledown$} at 5.4 1.5
\circulararc 75 degrees from 3.5  3  center at 3.5 1
\circulararc -75 degrees from 3.5  0  center at 3.5 2
\put{$\blacktriangle$} at 1.6 1.5
\circulararc -75 degrees from 3.5  3  center at 1.5 1.5
\circulararc 75 degrees from 3.5  3  center at 5.5 1.5
\put{$\blacktriangle$} at 4 1.5
\put{$\blacktriangledown$} at 3 1.5
\plot  -6 4   -1 4 /
\plot  -6 -1    -1 -1 /
\plot  -6 4   -6 -1 /
\plot  -1 4  -1  -1 /
\plot  1 4   6 4 /
\plot  1 -1    6 -1 /
\plot  1 4   1 -1 /
\plot  6 4  6  -1 /
%
%
\endpicture
\end{equation*}

The chambers $\widetilde{C}_i$ for $i \in \Z$ are the open intervals $(i,i+1)$ and the
$1$-codimensional facets are the points. The toric arrangement depends on the chosen lattice. For example we can quotient in two different way as in the above figure. 

Namely, the picture on the left corresponds to the choice $\Lambda=\Z$, i.e. $\pi:x\mapsto e^{2\pi i x}$, whereas the picture on the right is given by $\Lambda=2\Z$ and $\pi:x\mapsto e^{\pi i x}$.
 {As shown in the pictures the complex in the former example cannot be described by the two pairs
$[C_{-1} \prec C_{-1}]$, $[C_{-1} \prec e_0]$ since it has $3$ cells. 
Furthermore, this CW-complex is not regular (the closure of its cells is not contractible). On the other hand, in the latter example we have a regular CW-complex with two 0-dimensional cells and four 1-dimensional cells.

}

\bigskip

%

Now we will focus on the case in wich $\Cal S$ maps bijectively on the set of pairs 
$[C \prec F]$, since then the description of the complex $\mathcal{S}$ is particularly striking.
Since $\Cal S=\pi(\wdt{\Cal S})$ is a complex homotopic to the complement $\mathcal{R}_X$,
$\Cal S$ is described by  pairs of the form $[C \prec F]$  if and
only if the map
\begin{equation}\label{cond0}
\pi([\wdt C \prec \wdt F]) \longrightarrow  [\pi(\wdt C) \prec
\pi(\wdt F)] 
\end{equation}
is injective. 

\bigskip

Moreover, if the definition (\ref{cond0}) holds then we can define the boundary of a pair
$[C \prec F]$. We need first to introduce new notations.

 \bigskip


\textbf{Notations.} Let $P_0\subset V$ be a fundamental parallelogram for $\pi: V\rightarrow T$ containing the origin of $V$.
Let $\Cal{A}_{0,X}$ be the subarrangement of $\Cal{A}_{X}$ made by all the hyperplanes that intersect $P_0$ (see, for istance, figure (\ref{figura1}) in the next Section).

We will say that a maximal dimensional cell $[\wdt C \prec \wdt F^n]$ is in $\Cal{A}_{0,X}$ if its support $\mid \wdt F^n \mid$ is the intersection
of some of the hyperplanes in $\Cal{A}_{0,X}$. While a $k$-cell $[\wdt
C \prec \wdt F^k]$ is in $\Cal{A}_{0,X}$
if it is in the boundary of a $n$-cell in $\Cal{A}_{0,X}$. Let $\wdt {\Cal S}_{0}$ be the set of all
such cells.

\bigskip

With previous notations  if (\ref{cond0}) is injective (i.e. it is a bijection) we define the boundary 
as follow:

\bigskip

$[C \prec F^k]$ is in the boundary of $[D \prec G^j]$ $(k < j)$ if and only if there are cells
$[\wdt C \prec \wdt F^k] \in \pi^{-1}([C \prec F^k]) \cap \wdt {\Cal S}_{0}$
and $[\wdt D \prec \wdt G^j] \in \pi^{-1}([D \prec G^j]) \cap \wdt {\Cal S}_{0}$
such that $[\wdt C \prec \wdt F^k] \in \partial_{\wdt{\Cal S}} [\wdt D \prec \wdt G^j]$.

\bigskip

Obviously this boundary map commutes with the one in $\wdt{\Cal S}$ and
we get that the map in (\ref{cond0}) is a bijection of CW-complexes.

\bigskip

Toric arrangement for which $\Cal S$ is in bijection with pairs $[C \prec F]$ are easily characterized as follows.

\begin{df} A toric arrangement $\Cal{T}_X$ is \emph{thick} if the quotient map
$$\pi : V \longrightarrow T$$ is injective on the closure $clos(\widetilde{C})$ of every chamber
$\widetilde{C}$  of the associated affine arrangement $\Cal{A}_{X, \R}$.
\end{df}
We notice that every toric arrangement is covered by a thick one and the fiber of the covering map is finite;
hence our assumption is not very restrictive.

We have the following

\begin{lem}\label{isom} A toric arrangement $\Cal{T}_X$ is thick if and only if
$$
[\pi(\wdt C) \prec \pi(\wdt F)]=[\pi(\wdt D) \prec \pi(\wdt G)]  \Longleftrightarrow \pi([\wdt C \prec \wdt F])=\pi([\wdt D \prec \wdt G])
$$
for any two cells $[\wdt C \prec \wdt F], [\wdt D \prec \wdt G] \in \wdt{\Cal S}$
\end{lem}

\textbf{Proof.} By previous considerations, it is enough to prove that the thick condition is equivalent to
$$[\pi(\wdt C) \prec \pi(\wdt F)]=[\pi(\wdt D) \prec \pi(\wdt G)]  \Longrightarrow \pi([\wdt C \prec \wdt F])=\pi([\wdt D \prec \wdt G])$$
$\Rightarrow :$ Let $\Cal{T}_X$ be thick and $[\pi(\wdt C) \prec \pi(\wdt F)]=[\pi(\wdt D) \prec \pi(\wdt G)]$ for two given $k$-cells in $\wdt{\Cal S}$. This implies that $\pi(\wdt C) = \pi(\wdt D)$ and
$\pi(\wdt F)=\pi(\wdt G)$, i.e. there are translations $t,t^{\prime} \in \Lambda$ such that
$\wdt D=t.\wdt C$ and $\wdt G=t^{\prime}.\wdt F$.

By construction $t.\wdt F$ is a facet in the closure
$clos(D)$. We get two facets $t. \wdt F$ and $\wdt G$ both in $clos(D)$ and with the same image
$\pi(t.\wdt F)=\pi(\wdt F)=\pi(\wdt G)$. By hypothesis $\pi$ is injective on $clos(D)$ then
$t.\wdt F=G$, i.e. $t=t^{\prime}$ which implies that $\pi([\wdt C \prec \wdt F])=\pi([\wdt D \prec \wdt G])$.

$\Leftarrow$ Let $\wdt F$ and $\wdt G$ two facets in $clos(\wdt C)$ such that
$\pi(\wdt F)=\pi(\wdt G)$ then
$$\pi([\wdt C \prec \wdt F])=[\pi(\wdt C) \prec \pi(\wdt F)]=[\pi(\wdt C) \prec \pi(\wdt G)]=
\pi([\wdt C \prec \wdt G]).$$
As a consequence if $t \in \Lambda$ is the translation such that $\wdt F=t.\wdt G$ then
$t. \wdt C=\wdt C$. It follows that $t$ is the identity and we get $\wdt F=\wdt G$, i.e.
$\pi$ is injective on $clos(\wdt C)$ $\qquad$ $\square$

\bigskip

By Lemma \ref{isom} the map defined in (\ref{cond0}) is a bijection if and only if
$\mathcal{T}_X$  is a thick toric arrangement. Hence the set of pairs $[C \prec F]$ is a CW-complex $\overline{\Cal S}$ and we get the following theorem

\begin{teo}
Let  $\mathcal{T}_X$ be a thick toric arrangement.
Then its complement $\Cal{R}_X$ has the same
homotopy type of the CW-complex $\overline{\Cal S}$.
\end{teo}

Then in this case the complex $\mathcal{S}$ has a nice combinatorial description, totally analogue to that of the classical Salvetti complex \cite{Sa1}.

Moreover if a toric arrangement is thick then the maximal dimensional cells
$[\widetilde{C} \prec \widetilde{F}^n]$ in $\Cal{A}_{0,X}$ are in one to one correspondence with the $n$-dimensional facets of $\overline{\Cal{S}}$.
Then the boundary in a thick toric arrangement $\Cal T_X$ can be completely described knowing the boundary in the associated finite complex $\Cal A_{0,X}$.

This allows to better understand the fundamental group of the complement
and to perform computations on integer cohomology.

Furthermore, in this case $\mathcal{S}$ is a \emph{regular} CW-complex.

\begin{rmk} The number of chambers of $\Cal T_{X,\R}$ can be computed
by formulae given in \cite{ERS} and \cite{Mo2}. However, the
combinatorics of the layers in $\Cal T_{X,\R}$ is more complicated
than the one of spaces of  $\Cal A_{X,\R}$. Hence, an enumeration of
the facets is not easy to provide in the general case. Thus from now on
we focus on the arrangements defined by roots systems. In this case
the chambers are parametrized by the elements of the Weyl group, and
the poset of layers has been described in \cite{Mo1}.
\end{rmk}

\section{Weyl toric arrangements}

In this section we give a simpler description of the above complex for the case of toric arrangements associated to affine Weyl groups, by taking as $\Lambda$ the coroot lattice (for the theory of Weyl groups see, for instance, \cite{Bou}).
Indeed in this case the toric arrangement is thick.
Using this description, we give an example of how the integer cohomology of these arrangements can be computed.

\subsection{Notations and Recalls.}

\paragraph{Toric arrangement associated to a Weyl group.}
Let $\Phi$ be a root system, $\Lambda=\langle\Phi^\vee\rangle$ be the
lattice spanned by the coroots, and $\Delta$ be its dual
lattice (which is called the \emph{cocharacters} lattice).
Then we define a torus $T=T_\Delta$ having $\Delta$
as group of characters. Namely, if $\mathfrak{g}$
is the semisimple complex Lie algebra associated to $\Phi$ and
$\mathfrak{h}$ is a Cartan subalgebra, $T$ is defined as the
quotient $T\doteq {\mathfrak{h}}/\Lambda$.

Each root $\alpha$ takes integer values on
$\Delta$, so it induces a map
$$e^{\alpha}: T\rightarrow {\mathbb{C}}/{\mathbb{Z}}\simeq \mathbb{C^*}$$
which is a character of the torus.
Let $X$ be the set of these
characters; more precisely, since $\alpha$ and $-\alpha$ define
the same hypersurface, we set
$$X\doteq \left\{e^{\alpha},\:\alpha\in \Phi^+\right\}.$$
In this way to every affine Weyl group $\widetilde{W}$ we associate a toric
arrangement $\mathcal{T}_{\widetilde{W}}$, with complement $\mathcal{R}_{\widetilde{W}}$.

We will call these arrangements Weyl toric arrangements. They have been studied in \cite{Mcm} and \cite{Mo1}.

\begin{rmk}\label{remsec}

~

\begin{enumerate}
  \item Let $G$ be the semisimple, simply connected linear
      algebraic group associated to $\mathfrak{g}$. Then
      $T$ is the maximal torus of $G$ corresponding to
      $\mathfrak{h}$, and $\rx$ is known as the set of
      \emph{regular points} of $T$.
  \item One may take as $\Delta$ the root lattice (or equivalently, take as $\Lambda$ the character lattice).
 But in this way one obtains as $T$ a maximal
      torus of the semisimple \emph{adjoint} group $G^a$,
      which is the quotient of $G$ by its center.
\end{enumerate}
\end{rmk}

Let $(\widetilde{W},S)$ be the Coxeter system associated to $\widetilde{W}$ and
$$\mathcal{A}_{\widetilde{W}}=\{H_{\widetilde{w} s_i \widetilde{w}^{-1}} \mid \widetilde{w}
\in \widetilde{W} \mbox{ and } s_i \in S\}$$
the arrangement in $\C^n$ obtained by complexifying the reflection hyperplanes of
$\widetilde{W}$, where, in a standard way, the hyperplane
$H_{\widetilde{w} s_i \widetilde{w}^{-1}}$ is the hyperplane fixed by the reflection
$\widetilde{w} s_i \widetilde{w}^{-1}$.\\
We can view $\Lambda$ as a subgroup of $\widetilde{W}$, acting by translations. Then it is well known that $\widetilde{W}/\Lambda \simeq W$, where $W$ is the finite reflection group associated to $\widetilde{W}$ (see for instance \cite{Ra}). As a consequence, the toric arrangement
can be described as:
$$
T_{\widetilde{W}}=\{H_{[w] s_i [w^{-1}]} \mid w \in W \mbox{ and } s_i \in S \}
$$
where two hypersurfaces $H_{[w] s_i [w^{-1}]}$ and $H_{[\overline{w}] s_i [\overline{w}^{-1}]}$ are equal if and only if
there is a translation $t \in \Lambda$ such that $tw s_i (tw)^{-1}=\overline{w} s_i \overline{w}^{-1}$, i.e. $\overline{w}=tw$.

By \cite{Mo1}, these hypersurfaces
intersect in
$$
\frac{\mid W \mid}{\mid W_{S \setminus \{s_i \} } \mid}
$$
local copies of the finite hyperplane arrangement $A_{W_{S \setminus \{s_i\}}}$ associated to the group
 generated by $S \setminus \{s_i\}$,  $s_i \in S$.

For example in the affine Weyl group $\widetilde{A}_n$ generated by $\{s_0,\ldots, s_n\}$
for any generator $s_i$ the finite reflection group associated to $S \setminus \{s_i\}$ is a copy of the finite Coxeter group $A_n$.

Then we have
 {\begin{prop}The toric arrangement $\mathcal{T}_{\widetilde{W}}$ is thick. \end{prop} }

 {\textbf{Proof.} Since $\Lambda$ is the coroot lattice, if $t \in \Lambda$  is a translation such that there is a $n$-codimensional facet $\wdt F^n \in clos(\wdt C) \cap clos(t.\wdt C)$ for  an affine chamber $\wdt C$, then $t$ is the identity (see \cite{Bou}). }


 {If $T_{\wdt{W}}$ is not thick then there are two facets $\wdt F_1$ and $\wdt F_2$ in
the closure $clos(\wdt C)$ of a chamber $\wdt C$ such that $\pi(\wdt F_1)=\pi(\wdt F_2)$, i.e.
there is a translation $t \in \Lambda$ such that $\wdt F_2= t.\wdt F_1$. Hence $\wdt F_2$ is a facet
in $clos(C) \cap clos(t.C)$. In particular all the $n$-codimensional facets $\wdt F^n$ in the closure 
of $\wdt F_2$ are in the closure of both $C$ and $t.C$. This is a contradiction and it concludes the proof. $\qquad \square$}

\bigskip

Then we can construct
the Salvetti complex for these arrangements in a way which is very similar to the one known for affine Coxeter arrangements.

\paragraph{Salvetti Complex for affine Artin groups}
It is well known (see, for instance, \cite{boss3}, \cite{Sa} ) that the cells of Salvetti complex
$\widetilde {\Cal S}_W$ for arrangements $\Cal{A}_{\widetilde{W}}$ are of the form
$E(\widetilde{w},\Gamma)$ with
$\Gamma \subset S$ and $\widetilde{w} \in  \widetilde{W}$. Indeed if $\widetilde{\alpha} \in \{\widetilde{w}s\widetilde{w}^{-1} | s \in S, \widetilde{w} \in \widetilde{W}\}$ is a reflection,
the chambers are in one to one correspondence with the elements of the group $\widetilde{W}$ as follows.
Fixed a base chamber $C_0$, it corresponds to $1 \in \widetilde{W}$. Now if $C$ corresponds to $\widetilde{w}$, the chamber $D$ separated from $C$ by the reflection hyperplane $H_{\widetilde{\alpha}}$ corresponds to the element
$\widetilde{\alpha}\widetilde{w} \in \widetilde{W}$. The notation
$D \simeq \widetilde{\alpha}\widetilde{w}$ will be used.

If $\wdt F^k$ is a $k$-codimensional facet then the $k$-cell $[\wdt C \prec \wdt F^k]$
corresponds to the  pair $E(\widetilde{w},\Gamma)$ where $\widetilde{w} \simeq \wdt C$ and $\Gamma=\{s_{i_1},\ldots , s_{i_k}\}$ is the unique subset of cardinality $k$ in $S$ such that
$$\mid F^k \mid = \bigcap_{j=1}^k H_{\widetilde{w}s_{i_j}\widetilde{w}^{-1}}.$$
If $\widetilde{W}_{\Gamma}$ is the finite subgroup generated by $s \in \Gamma$, by \cite{boss3} the integer boundary map can be expressed as follows:
\begin{equation}
\begin{split}
\partial_k(E(\widetilde{w},\Gamma)) =
&\sum_{s_j \in \Gamma}
\sum_{\beta\in \widetilde{W}^{\Gamma\setminus\{ s_j \}}_{\Gamma}}(-1)^{l(\beta)+\mu
(\Gamma,s_j)} E(\widetilde{w}\beta,\Gamma\setminus\{ s_j\}).
\end{split}
\end{equation}
where $\widetilde{W}^{\Gamma\setminus\{\sigma\}}_{\Gamma}=\{w \in \wdt W_{\Gamma} : l(ws) > l(w) \forall s \in \Gamma \setminus\{\sigma\} \}$ and $\mu(\Gamma, s_j)=\sharp\{s_i \in \Gamma | i \leq j \}$.

\begin{rmk}Instead of the co-boundary operator we prefer to describe its dual,
i.e. we define the boundary of a $k$-cell $E(\widetilde{w},\Gamma )$ as a linear combination of
the $(k-1)$-cells which have $E(\widetilde{w},\Gamma )$ in their co-boundary, with the
same coefficient of the co-boundary operator.
We make this choice since the boundary operator has a nicer description than
co-boundary operator in terms of the elements of $\wdt W$.
\end{rmk}

\subsection{Description of the complex}

Let $\Cal{S}_W$ be the CW-complex associated to $\mathcal{T}_{\widetilde{W}}$. By the previous considerations, $\Cal{S}_W$ admits a description similar to that of
$\widetilde{\Cal S}_W$.
Indeed each chamber $C$ is in one to one correspondence with an equivalence class
$[w] \in \widetilde{W} / \Lambda$ and then with an element
$w \in W \simeq \widetilde{W} / \Lambda$ of the finite reflection group $W$. We will write $C \simeq [w]$.

In the same way, the  pair $[C \prec F^k]$
corresponds to the cell $E([w],\Gamma) \in \Cal{S}_W$ where  $C \simeq [w]$ and
$\Gamma=\{s_{i_1},\ldots , s_{i_k}\}$ is the unique subset of cardinality $k$ in $S$ such that
$$\mid F^k \mid = \bigcap_{j=1}^k H_{[w]s_{i_j}[w^{-1}]}.$$

We now want to describe the boundary of each cell: this is done in a
standard way by characterizing the cells that are in the boundary of a
given cell, and by assigning an orientation to all cells (see, for
instance, \cite{Sa}).

By construction the toric CW-complex is locally isomorphic to the affine one and it can inherit its affine orientation. Then the integer boundary operator for Weyl toric arrangements can be written as the affine one:
\begin{equation}\label{bordo}
\begin{split}
\partial_k(E([w],\Gamma)) =
&\sum_{\sigma\in \Gamma}
\sum_{\beta\in W^{\Gamma\setminus\{\sigma\}}_{\Gamma}}(-1)^{l(\beta)+\mu
(\Gamma,\sigma)} E([w\beta],\Gamma\setminus\{ \sigma \})
\end{split}
\end{equation}
where, instead of elements of the affine group $\widetilde{W}$, we have equivalence classes with representatives in the finite group $W$. 

By the formula above, the complex $\Cal{S}_W$ can be effectively used for computing homotopy invariants of $\mathcal{R}_{\widetilde{W}}$. For instance we have

 {\begin{prop} $$H^{\bullet}(\mathcal{R}_{\widetilde{W}},\Z) \simeq H^{\bullet}(\Cal{S}_W,\Z) $$
where the coboundary map is the dual of the map defined in (\ref{bordo}).
\end{prop}}


\bigskip

\textbf{Example.} Let us consider the affine Weyl group $\widetilde{B}_2$ (see \cite{Bou}) with Coxeter-Dynkin diagram

\begin{equation*}
\begin{array}{ccccc}
\circ & \stackrel{4}{-} &\circ & \stackrel{4}{-} &\circ \\
s_0& & s_1& &s_2 \\
\end{array}
\end{equation*}

and associated finite group $B_2$
\begin{equation*}
\begin{array}{ccc}
\circ & \stackrel{4}{->} &\circ \\
s_1& &s_2 \\
\end{array}
\end{equation*}

In this case we get translations $t_1=s_0s_1s_2s_1$ and $t_2=s_2s_1s_0s_1$ and the affine arrangement is represented as:
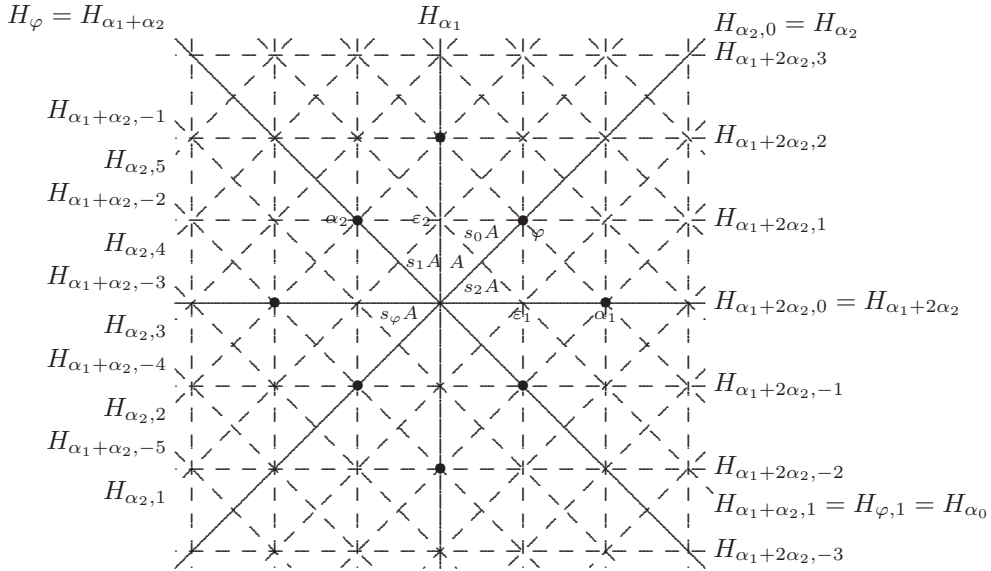
\begin{equation*}
\beginpicture
\setcoordinatesystem units <1.1cm,1.1cm>         
\setplotarea x from -4.5 to 4.5, y from -3.5 to 4    
\put{$H_{\alpha_1}$}[b] at 0 3.3
\put{$H_{\alpha_2,0}=H_{\alpha_2}$}[bl] at 3.3 3.2
\put{$H_{\alpha_2,1}$}[r] at -3.3 -2.3
\put{$H_{\alpha_2,2}$}[r] at -3.3 -1.3
\put{$H_{\alpha_2,3}$}[r] at -3.3 -0.3
\put{$H_{\alpha_2,4}$}[r] at -3.3 0.7
\put{$H_{\alpha_2,5}$}[r] at -3.3 1.7
\put{$H_\varphi = H_{\alpha_1+\alpha_2}$}[br] at -3.3 3.3
\put{$H_{\alpha_1+\alpha_2,1}=H_{\varphi,1}=H_{\alpha_0}$}[tl] at 3.3 -2.3
\put{$H_{\alpha_1+\alpha_2,-1}$}[r] at -3.3 2.3
\put{$H_{\alpha_1+\alpha_2,-2}$}[r] at -3.3 1.3
\put{$H_{\alpha_1+\alpha_2,-3}$}[r] at -3.3 0.3
\put{$H_{\alpha_1+\alpha_2,-4}$}[r] at -3.3 -0.7
\put{$H_{\alpha_1+\alpha_2,-5}$}[r] at -3.3 -1.7
\put{$H_{\alpha_1+2\alpha_2,0}=H_{\alpha_1+2\alpha_2}$}[l] at 3.3 0
\put{$H_{\alpha_1+2\alpha_2,1}$}[l] at 3.3 1
\put{$H_{\alpha_1+2\alpha_2,2}$}[l] at 3.3 2
\put{$H_{\alpha_1+2\alpha_2,3}$}[l] at 3.3 3
\put{$H_{\alpha_1+2\alpha_2,-1}$}[l] at 3.3 -1
\put{$H_{\alpha_1+2\alpha_2,-2}$}[l] at 3.3 -2
\put{$H_{\alpha_1+2\alpha_2,-3}$}[l] at 3.3 -3
\put{$\scriptstyle{A}$} at 0.2 0.5
\put{$\scriptstyle{s_1A}$} at -0.2 0.5
\put{$\scriptstyle{s_2A}$} at 0.5 0.2
\put{$\scriptstyle{s_0A}$} at 0.5 0.85
\put{$\scriptstyle{s_\varphi A}$} at -0.5 -0.15
\put{$\scriptstyle{\varepsilon_1}$}[t] at 1 -0.1
\put{$\scriptstyle{\varepsilon_2}$}[r] at -0.1 1
\put{$\scriptstyle{\alpha_1}$}[t] at   2  -0.1
\put{$\scriptstyle{\alpha_2}$}[r] at -1.1   1
\put{$\scriptstyle{\varphi}$}[tl] at  1.1  0.9
\put{$\bullet$} at  0  2
\put{$\bullet$} at -1  1
\put{$\bullet$} at  1  1
\put{$\bullet$} at -2  0
\put{$\bullet$} at  2  0
\put{$\bullet$} at -1 -1
\put{$\bullet$} at  1 -1
\put{$\bullet$} at  0 -2
\plot -3.2 -3.2   3.2 3.2 /
\plot  3.2 -3.2  -3.2 3.2 /
\plot  0  3.2   0 -3.2 /
\plot  3.2  0  -3.2  0 /
\setdashes
\plot  -3 3.2   -3 -3.2 /
\plot  -2 3.2   -2 -3.2 /
\plot  -1 3.2   -1 -3.2 /
\plot   1 3.2    1 -3.2 /
\plot   2 3.2    2 -3.2 /
\plot   3 3.2    3 -3.2 /
\plot  3.2 -3   -3.2 -3 /
\plot  3.2 -2   -3.2 -2 /
\plot  3.2 -1   -3.2 -1 /
\plot  3.2  1   -3.2  1 /
\plot  3.2  2   -3.2  2 /
\plot  3.2  3   -3.2  3 /
\plot  3.2 -1.8   1.8 -3.2 /
\plot  3.2 -0.8   0.8 -3.2 /
\plot  3.2 0.2    -0.2 -3.2 /
\plot  3.2 1.2   -1.2 -3.2 /
\plot  3.2 2.2   -2.2 -3.2 /
\plot  2.2 3.2   -3.2 -2.2 /
\plot  1.2 3.2   -3.2 -1.2 /
\plot  0.2 3.2   -3.2 -0.2 /
\plot  -0.8 3.2  -3.2 0.8 /
\plot  -1.8 3.2  -3.2 1.8 /
\plot  -2.2 3.2   3.3 -2.3 /
\plot  -1.2 3.2   3.2 -1.2 /
\plot  -0.2 3.2   3.2 -0.2 /
\plot  0.8 3.2  3.2 0.8 /
\plot  1.8 3.2  3.2 1.8 /
\plot  -3.2 -1.8   -1.8 -3.2 /
\plot  -3.2 -0.8   -0.8 -3.2 /
\plot  -3.2 0.2    0.2 -3.2 /
\plot  -3.2 1.2   1.2 -3.2 /
\plot  -3.2 2.2   2.2 -3.2 /
\endpicture
\end{equation*}
If $\Cal{A}_0$ is the finite subarrangement defined in Section 2.3, then the real toric arrangement is obtained quotienting it as shown in the
following figure, where arrows indicate identified edges:

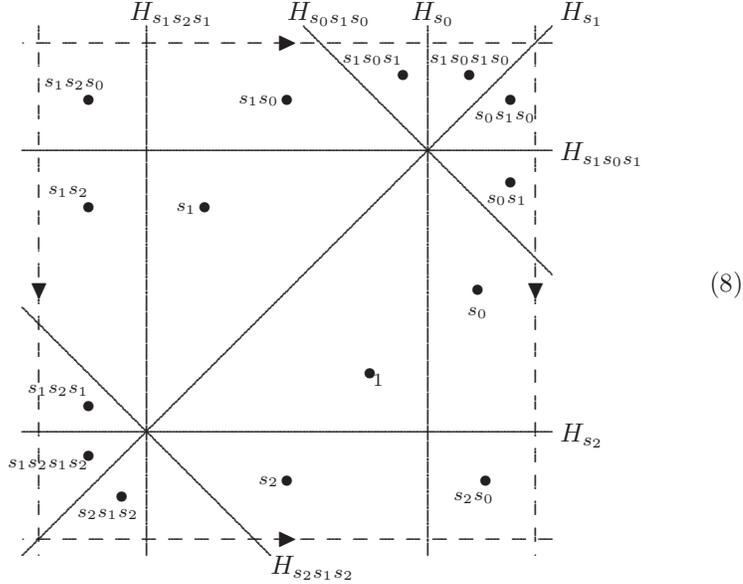
\begin{equation} \label{figura1}
\beginpicture
\setcoordinatesystem units <1.1cm,1.1cm>         
\setplotarea x from -4.5 to 4.5, y from -3.5 to 4    
\put{$H_{s_1}$}[bl] at 3.3 3.2
\put{$H_{s_0}$}[bl] at 1.5 3.2
\put{$H_{s_1s_2s_1}$}[bl] at -1.9 3.2
\put{$H_{s_0s_1s_0}$}[bl] at 0 3.2
\put{$H_{s_1s_0s_1}$}[bl] at 3.3 1.5
\put{$H_{s_2s_1s_2}$}[bl] at -0.2 -3.5
\put{$H_{s_2}$}[bl] at 3.3 -1.9
\put{$\bullet$} at 2.3 0
\put{$\scriptstyle{s_0}$}[t] at 2.3 -0.2
\put{$\bullet$} at 0 2.3
\put{$\scriptstyle{s_1s_0}$}[r] at -0.1 2.3
\put{$\bullet$} at  2.7  2.3
\put{$\scriptstyle{s_0s_1s_0}$}[r] at 3 2.1
\put{$\bullet$} at  2.2  2.6
\put{$\scriptstyle{s_1s_0s_1s_0}$}[r] at  2.7 2.8
\put{$\bullet$} at  1.4  2.6
\put{$\scriptstyle{s_1s_0s_1}$}[r] at 1.4 2.8
\put{$\bullet$} at 2.7 1.3
\put{$\scriptstyle{s_0s_1}$}[r] at 2.9 1.1
\put{$\bullet$} at  -2.4  2.3
\put{$\scriptstyle{s_1s_2s_0}$}[r] at -2.2 2.5
\put{$\bullet$} at -2.4 1
\put{$\scriptstyle{s_1s_2}$}[r] at -2.4 1.2
\put{$\bullet$} at -2.4 -1.4
\put{$\scriptstyle{s_1s_2s_1}$}[r] at -2.4 -1.2
\put{$\bullet$} at -2.4 -2
\put{$\scriptstyle{s_1s_2s_1s_2}$}[r] at -2.4 -2.1
\put{$\bullet$} at -2 -2.5
\put{$\scriptstyle{s_2s_1s_2}$}[r] at -1.8 -2.7
\put{$\bullet$} at 1 -1
\put{$\scriptstyle{1}$}[t] at  1.1 -1
\put{$\bullet$} at -1 1
\put{$\scriptstyle{s_1}$}[r] at -1.1 1
\put{$\bullet$} at  2.4  -2.3
\put{$\scriptstyle{s_2s_0}$}[r] at 2.5 -2.5
\put{$\bullet$} at  0 -2.3
\put{$\scriptstyle{s_2}$}[r] at -0.1 -2.3
\plot -3.2 -3.2   3.2 3.2 /
\plot  -1.7 3.2   -1.7 -3.2 /
\plot   1.7 3.2    1.7 -3.2 /
\plot  3.2 -1.7   -3.2 -1.7 /
\plot  3.2  1.7   -3.2  1.7 /
\plot  0.2 3.2  3.2 0.2 /
\plot  -3.2 -0.2   -0.2 -3.2 /
\put{$\blacktriangleright$} at 0 3
\put{$\blacktriangleright$} at 0 -3
\put{$\blacktriangledown$} at -3 0
\put{$\blacktriangledown$} at 3 0
\setdashes
\plot  -3 3.2   -3 -3.2 /
\plot   3 3.2    3 -3.2 /
\plot  3.2 -3   -3.2 -3 /
\plot  3.2  3   -3.2  3 /
\endpicture
\end{equation}
Here, for brevity,  the vertices $E(w, \emptyset)$ are labelled by the element
 $w \in \wdt W$.

 \bigskip

We get, for example, that the cell $E([1],\emptyset)$ is the vertex in the chamber containg
$1 \in \wdt W$, while the vertices $E([s_0],\emptyset)$ and
 $E([s_1s_2s_1],\emptyset)$ correspond to the same chamber in the toric arrangement; indeed
 $s_0=t_1s_1s_2s_1$, then $[s_0]=[s_1s_2s_1]$.

Notice that the number of chambers in the real torus is $8$ in one to one correspondence with the finite Weyl group $B_2$ with cardinality $8$.
Then we get exactly:

\bigskip

 $8 \qquad 0$-cells of the form $E([w],\emptyset)$ for $w \in B_2$,

 $24 \qquad 1$-cells of the form $E([w],\{s_i\})$ for $w \in B_2$ and $i=0,1,2$,

 $24 \qquad 2$-cells of the form $E([w],\{s_i,s_j\})$ for $w \in B_2$ and $0 \leq i<j \leq 2$.

 \bigskip

These cells locally correspond to four finite Coxeter arrangements, two of type $B_2$ and two of type $A_1 \times A_1$
appearing in the figure above.
 In particular  the $2$-cells can be written as:

 \bigskip

  $E([w],\{s_i,s_{i+1}\})$ with a representative $w$ chosen in the Coxeter group $B_2$ generated by $\{s_i,s_{i+1}\})$, i=0,1;

 $E([w],\{s_0,s_2\})$ and $E([s_1w],\{s_0,s_2\})$ with a representative $w$ chosen in the group
  $\{1,s_0,s_2,s_0s_2\}$ generated by $\{s_0,s_2\}$.

\bigskip

 The representatives can be chosen in the more suitable way for computations. The boundary
 map (\ref{bordo}) for the $1$-cells is:
\begin{equation*}
\partial_1 E([w],\{s_i\})= E([w],\emptyset) - E([ws_i],\emptyset)
\end{equation*}
and it gives rise to a matrix of $24$ columns and $8$ rows with entries $0$, $1$ and $-1$.\\
On the other hand, the second boundary map is given by
 \begin{equation*}
 \begin{split}
 \partial_2 E([w],\{s_i,s_{i+1}\})= E([w],\{s_i\}) - E([ws_{i+1}],\{s_i\}) + E([ws_{i}s_{i+1}],\{s_i\})- \\
 - E([w],\{s_{i+1}\})  + E([ws_{i}],\{s_{i+1}\}) - E([ws_{i+1}s_{i}],\{s_{i+1}\})\\
 \end{split}
 \end{equation*}
  \begin{equation*}
  \partial_2 E([w],\{s_0,s_2\})= E([w],\{s_0\}) - E([ws_2],\{s_0\}) - E([w],\{s_{2}\})  + E([ws_{0}],\{s_2\}).
 \end{equation*}
In this way we get that the homology, and hence the cohomology, is torsion free and $H_0(R_{B_2}, \Z) = \Z$, $H_1(R_{B_2}, \Z) = \Z^8$ and $H_2(R_{B_2}, \Z) = \Z^{15}$, which agrees with the Betti numbers computed in \cite[Ex. 5.14]{Mo1}.

In general we have the following

\begin{conj}
Let $\widetilde{W}$ be an affine Weyl group and $\mathcal{T}_{\widetilde{W}}$ be the corresponding toric arrangement. Then the integer cohomology of the complement is torsion free (and hence it coincides with the De Rham cohomology computed in \cite{DP}).
\end{conj}

This conjecture will be proved in a future paper \cite{simo2}.

\section{An example from robotics}

In this section we give an example of non-thick arrangement: the one obtained from the
affine Weyl arrangement $\Cal A_{\wdt A_n}$, by quotienting by the coroot lattice, which we will denoted by
$\Lambda_{\wdt A_n}$ (see the second part of Remark \ref{remsec}).

\bigskip

Indeed in this case the underlying real toric arrangement has a very nice description in terms of Young tableaux. More precisely the facets of $\Cal T_{\wdt A_n, \R}$ are in one to one correspondence with
a family of Young tableaux which turn out to be the same tableaux describing the space of all periodic legged gaits of a robot body (see \cite{KGCohen}).

\bigskip

It is clear that, in this case, the finite arrangement $\Cal A_{0,\wdt A_n}$ is exactly the braid
arrangement $\Cal A_{A_n}$.

\subsection{Tableaux description for the complex $\wdt{\Cal S}_{A_n}$}

We indicate by $A_n$ the symmetric group on $n+1$ elements,
acting by permutations of the coordinates. Then $\Cal{A} =
\Cal{A}_{A_n}$ is the braid arrangement and $\wdt{\Cal S}_{A_n}$ is the
associated CW-complex (even if the arrangement is finite we continue to use the same notation
used above for the affine case to distinguish it from the toric one).

Given a system of coordinates in $\R^{n+1}$, we describe
$\wdt{\Cal S}_{A_n}$ through certain tableaux as follow.

Every $k$-cell $[\wdt C \prec \wdt F]$ is represented by a tableau with $n+1$
boxes and  $n+1-k$ rows (aligned on the left), filled with all the
integers in $\{1,...,n+1\}.$ There is no monotony condition on the
lengths of the rows. One has:
\medskip

\ni - $(x_1,\ldots, x_{n+1})$ is a point in $\wdt F$ if and only if:

\bigskip

$1.$ $i$ and $j$ belong to the same row if and only if $x_i=x_j$,

$2.$ $i$ belongs to a row preceding the one containing $j$ if and only if
$x_i < x_j$;

\bigskip

\ni - the chamber $\wdt C$ belongs to the half-space $x_i < x_j$ if and only if:

\bigskip

$1.$ either the row which contains $i$ is preceding the one
containing $j$ or

$2.$ $i$ and $j$ belong to the same row and the column which
contains $i$ is preceding the one containing $j$.
\medskip

Notice that the facets of the real stratification are represented
by standard Young tableaux, since the order of the entries in each
row does not matter, and hence we can assume it to be strictly increasing.\\
Notice also that the geometrical action of $A_n$ on the stratification
induces a natural action on the complex $\wdt{\Cal S}_{A_n}$ which, in terms of
tableaux, is given by a left action of $A_n$: $\sigma. \ T$ is the
tableau with the same shape as $T,$ and with entries permuted
by $\sigma.$

\subsection{Tableaux description for the facets of $\Cal T_{\wdt{A}_n, \R}$}

Let $\Cal A_{0,\wdt A_n} \subset \Cal A_{\wdt A_n}$ be the braid arrangement passing through the origin and $\pi: \R^{n+1} \longrightarrow \R^{n+1}/\Lambda_{\wdt A_n}=
T_{\R}$ the projection map.

\bigskip

If $\mathbf F_{\wdt A_n}$ is the stratification of $\R^{n+1}$ into facets induced by the arrangement
$\Cal A_{\wdt A_n}$, we define the set:
$$
\mathbf F_{0, \wdt A_n}=\{\wdt F^k \in \mathbf F_{\wdt A_n} \mid clos(F^k) \supset
\bigcap_{H \in \Cal A_{0,\wdt A_n} } H\}.
$$
Obviously $\mathbf F_{0,\wdt A_n}$ is in one to one correspondence with the stratification
$\mathbf F_{A_n}$ induced by the braid arrangement $\Cal A_{A_n}$ and the restriction
$\pi_{\mathbf F_{0,\wdt A_n}}$ is surjective on $T_{\R}$.

\bigskip

It follows that in order to understand how $\Lambda_{\wdt A_n}$ acts on
$\mathbf F_{\wdt A_n}$ it is enough to study how it acts on $\mathbf F_{0,\wdt A_n}$.
Moreover it is enough to consider facets in the closure of the base chamber $\wdt C_0$
corresponding to $1 \in \wdt A_n$; the action on the others will be obtained by symmetry.

\bigskip

Let us remark that a facet $\wdt F^k$ is in $\mathbf F_{0,\wdt A_n}$ if and only if it intersects any ball $B_0$
around the origin. Let $B_0$ be a ball of sufficiently small radius and
$$x=(x_1, \ldots, x_{n+1}) \in clos(\wdt C_0) \cap B_0$$ be a given point in a facet
$\wdt F^k \in \mathbf F_{0,\wdt A_n}$.
 Then the $x_i$'s satisfy $x_1 \leq x_2 \leq \ldots \leq x_{n+1}$ and the standard Young
tableaux $Tb_{\wdt F^k}$ associated to $\wdt F^k$ will have entries increasing along both,
rows and columns.

\bigskip

Let $t_1, \ldots ,t_n \in \Lambda_{\wdt A_n}$ be a base such that $t_i$ translates the
reflection hyperplane $H_{i,i+1}=Ker(x_i - x_{i+1})$ fixing all hyperplanes
$H_{j,j+1}=Ker(x_j-x_{j+1})$ for $j \neq i$  (i.e. each point in $H_{j,j+1}$ is sent in a point still in $H_{j,j+1}$).\\
Then we can assume that translation $t_i$ acts on the entry $x_i$ as $t_i.x_i= x_i + t$ with $x_i +t > x_{i+1}$ and, as $H_{j,j+1}$, for $j \neq i$, are invariant under the action of $t_i$, it follows that
$t_i.x_{i-1}=x_{i-1}+t$ and, by induction, $t_i.x_j=x_j+t$ for all $j < i$, while
$t_i.x_j=x_j$ for all $j >i$.

\bigskip

Recall that, by construction, given a standard Young tableaux, a point $(x_1,\ldots, x_{n+1})$ is a point in $\wdt F$ if and only if:

\bigskip

$1.$ $i$ and $j$ belong to the same row if and only if $x_i=x_j$,

$2.$ $i$ belongs to a row preceding the one containing $j$ if and only if $x_i < x_j$;

\bigskip

It follows that if $Tb$ is a tableau such that $i \in r_k$ and
$i+1 \in r_{k+1}$ are in two different rows, then $t_i$ acts on $Tb$ sending it in a tableau
$Tb^{\prime}$ with rows $r_1^{\prime}=r_{k+1}, \ldots, r_{h-k}^{\prime}=r_h,r_{h-k+1}^{\prime}=r_1, \ldots, r_h^{\prime}=r_k$. While if $i, i+1 \in r_k$ are in the same row, then $t_i$ acts sending the corresponding facet in a facet which is not anymore in $\Cal A_{0, \wdt A_n}$.

\bigskip

Then $\Lambda_{\wdt A_n}$ acts on the $h$ rows of a tableau $Tb_{\wdt F}$ as a power of the cyclic permutation $(1, \ldots , h)$.

\bigskip

Equivalently let $Y(n+1,k+1)$ be the set of standard Young tableaux with $k+1$ rows and $n+1$ entries and $Tb \in Y(n+1,k+1)$ be a tableau of rows $(r_1, \ldots , r_{k+1})$.
Then  {we have the following proposition.}

 {\begin{prop} The set of facets $F^k$ of the toric arrangement $\Cal T_{\wdt A_n,\R}$ is in one to one correspondence with the set
$$Y(n+1,k+1) / \sim$$
where a tableau $Tb^{\prime} \sim Tb$ if and only if the rows of $Tb^{\prime}$ are
$(r_{\sigma^s(1)}, \ldots , r_{\sigma^s(k+1)})$ for a power $\sigma^s$ of the cyclic permutation
$\sigma=(1, \ldots , k+1)$. 
\end{prop}}
In this way we get exactly the tableaux described in \cite{KGCohen}.

\bigskip

Finally let us recall that the relation $\wdt F^k \prec \wdt F^{k+1}$
holds if and only if the tableau $Tb_{\wdt F^{k+1}}$
corresponding to $\wdt F^{k+1}$ is obtained by attaching two consecutive rows of
$Tb_{\wdt F^k}$.\\
As a consequence if $F^k$ and $F^{k+1}$ are facets in the toric arrangement
$\Cal T_{\wdt A_n, \R}$,  $ F^k \prec F^{k+1}$ if and only if  the tableau $Tb_{F^{k+1}}$
corresponding to $F^{k+1}$ is obtained by attaching two consecutive rows of
$Tb_{F^k}$ or attaching the first one to the last one.

\end{document}